\documentclass[a4paper,11pt]{amsart}
\usepackage[utf8]{inputenc}
\usepackage{amsfonts,amssymb}
\usepackage{graphicx,tikz}
\usepackage{tikz-cd} 
\usepackage{float}
\usepackage{amsmath, amsthm, amsfonts}
\usepackage{geometry}
\usepackage{hyperref}
\usepackage{enumitem} 
\usepackage[capitalise]{cleveref}
\usepackage{comment}
\usepackage{amsrefs}
\usepackage[official]{eurosym}
\usepackage{graphicx}
\usepackage{nomencl}
\usepackage{amsfonts}
\usepackage{hyperref}
\usepackage{amssymb}
\usepackage{fancyhdr}
\usepackage{amscd}
\usepackage[english]{babel}

\theoremstyle{plain}
\theoremstyle{definition}

\newtheorem*{theorem*}{Theorem}
\newtheorem{thm}{Theorem}
\newtheorem{prop}[thm]{Proposition}

\theoremstyle{definition}
\newtheorem{defn}[thm]{\scshape{Definition}}

\DeclareMathOperator{\Aut}{Aut}
\DeclareMathOperator{\End}{End}

\DeclareMathOperator{\Hom}{Hom}
\DeclareMathOperator{\Prob}{Prob}
\makeatletter
\def\cleardoublepage{\clearpage\if@twoside \ifodd\c@page\else
	\hbox{}
	\thispagestyle{empty}
	\newpage
	\if@twocolumn\hbox{}\newpage\fi\fi\fi}
\makeatother

\title[Algorithmic problems in Engel groups]{Algorithmic problems in Engel groups and Cryptographic Applications}
\author[D. Kahrobaei]{Delaram Kahrobaei}
\address{University of York, Department of Computer Science, United Kingdom; The City University of New York, CUNY Graduate Center, U.S.A.; Computer Science and Engineering Department, New York University, U.S.A.}
\email{delaram.kahrobaei@york.ac.uk, dkahrobaei@gc.cuny.edu, dk2572@nyu.edu}
\author[M. Noce]{Marialaura Noce}
\address{Dipartimento di Matematica, University
of Salerno, Salerno, Italy; Matematika Saila, University of the Basque Country, Bilbao, Spain; Department of Mathematical Sciences, University of Bath, United Kingdom}
\email{mn670@bath.ac.uk}

\keywords{Engel elements, algorithmic problems, cryptography}
\subjclass[2010]{20F45, 94A60}

\begin{document}

\maketitle

\begin{abstract}
The theory of Engel groups plays an important role in group theory since they are closely related to the Burnside problems.

In this survey we consider several classical and novel algorithmic problems for Engel groups and propose several open problems. We study these problems with a view towards applications to cryptography.
\end{abstract}

\vskip 0.2 true cm

%------------------------------------------------------------------------------------%

%------------------------------------------------------------------------------------%
%------------------------------------------------------------------------------------%

\section{\bf Introduction}

In cryptography most common protocols (RSA, Diffie-Hellman, and elliptic curve methods) depend on the structure of commutative groups and they are related to the difficulty to solve integers factorization and discrete logarithms. In 1994 Shor provided a quantum algorithm that solves these problems in polynomial time \cite{Shor}. For this reason, researchers are motivated to find alternative methods for constructing cryptosystems. One of them is based on non-commutative cryptography, which does not operate over the integers. Hence, for security reasons, in the last decade new cryptosystems and key exchange protocols based on non-commutative cryptographic platforms have been developed.

The complexity of algorithmic problems have made available families of groups as platform groups for cryptographic protocols.  Among others, we mention braid groups (using the conjugacy search problem \cite{braid}), polycyclic groups (\cite{Eickconj} and \cite{polycyclic}), linear groups \cite{lineargroups}, and right-angled Artin groups (\cite{flores} and \cite{raags}).
For this reason, the employment of algorithmic group theoretic problems in cryptography is an active area of research nowadays.

In this paper we will present the actual state of Engel-group based cryptography. In particular, we investigate several group theoretic problems in Engel groups, with a view towards their applications to cryptography via computational complexity. We are primarily motivated by the fact that some of these group theoretic problems can be used for cryptographic purposes, such as authentication schemes, secret sharing schemes, key exchange problems, and multilinear maps.

The study of Engel groups originates from the famous paper of Burnside in 1901 and it is closely related to the General Burnside Problem (in the reminder GBP, for short)  \cite{Burnside}.
The General Burnside Problem asks if a finitely generated torsion group is finite. 
The question whether a finitely generated ($n$-)Engel group is nilpotent is the analogue of the GBP in the realm of Engel groups. In the following, we refer to this question as the ``Engel Burnside Problem''. Of course, every group that is locally nilpotent is an Engel group, but the converse needs not to be true in general. The first main result on Engel groups is the Theorem of Zorn which shows that any finite Engel group is nilpotent. 

We recall that an element $g$ of a given group $G$ is said to be \textit{right Engel} if for every $x \in G$ there exists an integer $n = n(g, x) \geq 1$ such that $[g, {}_nx] = 1$, where the commutator $[g, {}_nx]$ is defined recursively by the rules $[g, x] = g^{-1}g^{x}$ and 
$$
[g,_n x]=[g, x, \overset{n}{\dots}, x] = [[g, x, \overset{n-1}{\dots}, x], x]
$$
if $n > 1$. Similarly, $g$ is a \textit{left Engel} element if the variable $x$ appears on the left.
 A group is said to be an Engel group if all its elements are right Engel or, equivalently, left Engel.  
In a similar way, one can talk about \textit{bounded Engel elements}, with the condition that the integer $n$ can be chosen independently of $x$. A group $G$ is $n$-Engel if for any $x, y \in G$ we have $[x,_n y] =1$ for some $n \in \mathbb{N}$.

The paper is organized as follows. In Section \ref{Engelgroups} we give basic definitions and properties of Engel groups. We provide some identities satisfied by $n$-Engel groups, and, in particular, we survey known and open questions about the Engel Burnside Problem. We conclude the section by describing Burnside groups since we think that they have similar behaviour to Engel groups.

Section \ref{otherproperties} is devoted to survey general problems and interesting properties in the setting of Engel groups. We present the definition of degree of nilpotency for finite and infinite groups and some useful results that show how one can compute it. In addition, we give the definition of the growth rate and the Discrete Logarithm Problem.

In Section \ref{Sectionalgorithmicproblems} we introduce several  group-theoretic decision problems including the word, conjugacy, and isomorphism decision problem. We also describe the power decision problem, the endomorphism decision problem, the $n$-th root decision problem and the geodesic length decision problem. We present also many search problems.

In Section \ref{applications} we describe a number of cryptosystems that have been built around Engel groups. In particular, we present multilinear maps, a public key based on $2$-Engel groups and a digital signature based on $4$-Engel groups. In addition, we present two secret sharing schemes based on the efficiency of the word problem, and a key protocol whose platform groups are free nilpotent $p$-groups. Moreover, we survey the Learning Problem Homomorphism which uses as platform group the Burnside groups of exponent $3$ (that are, in particular, $2$-Engel groups). We conclude the section by presenting the Discrete Logarithm Problem in finite $p$-groups and the status of quantum algorithms.

The last Section \ref{conclusion} of the paper is devoted to survey the current status of Engel groups in cryptography. We also include a list of open problems, which we hope will guide researchers who wish to work in this  field.

\section{\bf Background on Engel groups}\label{Engelgroups}

\subsection{Preliminaries}\label{preliminaries}
Let $G$ be a group and let $x_{1}, x_{2}, \dots$ be elements of $G$. We define the commutator of weight $n \geq 1$ recursively by the rule
\begin{center}
	$[x_{1}, \dots, x_{n}] = \begin{cases} x_{1} & \mbox{if }n = 1 \\ [[x_{1}, \dots, x_{n-1}], x_{n}] & \mbox{if }n > 1,
	\end{cases}$
\end{center}
where $[x_{1}, x_{2}] = x_{1}^{-1}x_{2}^{-1}x_{1}x_{2} = x_{1}^{-1}x_{1}^{x_{2}}$.
If $x = x_{1}$ and $y = x_{2} = \dots = x_{n+1}$, we use the shorthand notation \[[x, {}_ny] = [x, y, \overset{n}{\dots}, y].\]

In the following we recall the definition of Engel elements of a group.
\begin{defn}
Let $G$ be a group and $g \in G$. We say that $g$ is a \textit{right Engel element} if for any $x \in G$ there exists $n = n(g, x) \geq 1$ such that $[g, {}_nx] = 1.$ If $n$ can be chosen independently of $x$, then $g$ is a \textit{right n-Engel element} (or a \textit{bounded right Engel element}).

Similarly, $g$ is a \textit{left Engel element} if for any $x \in G$ there exists $n = n(g, x) \geq 1$ such that $[x, {}_ng] = 1$. Again, if the choice of $n$ is independent of $x$, then $g$ is a \textit{left n-Engel element} (or a \textit{bounded left Engel element}).
\end{defn}

The sets of right and left Engel elements of $G$ are denoted by $R(G)$ and $L(G)$, respectively. Notice that $R(G)$ and $L(G)$ are invariant subsets under automorphisms of $G$. Hence, for any $g \in R(G)$ (resp. $g \in L(G)$) and any $x \in G$, we have $g^{x} \in R(G)$ (resp. $g^{x} \in L(G)$).

The following gives the iteration between right and left Engel elements.

\begin{prop}[\cite{hein}]\label{inclusion}
Let $G$ be a group. We have 
$$
R(G)^{-1} \subseteq L(G) \mbox{ and } R_{n}(G)^{-1} \subseteq L_{n+1}(G).
$$
\end{prop}

Notice also that in every $2$-group, elements of order $2$ are left Engel, as the following proposition (it can be found, for example, in \cite{rob}) shows.

\begin{prop}\label{inv}
	Let $G$ be  group and let $g \in G$ such that $g^{2} = 1$. Then for any $x \in G$ and any $n \geq 1$: \[[x, {}_ng] = [x, g]^{(-2)^{n-1}}.\] In particular, every involution in any $2$-group is a left Engel element.
\end{prop}

    In the following we give definitions of Engel and $n$-Engel groups.
    
\begin{defn}
We say that $G$ is an \textit{Engel group} if $G = R(G)$ or, equivalently, $G = L(G)$. Moreover, $G$ is an \textit{$n$-Engel} group if there exists $n \geq 1$ such that $[x, {}_ny] = 1$, for all $x, y \in G$.
\end{defn}

Of course every $n$-Engel group is Engel, and every nilpotent group of class $n$ is $n$-Engel. Also, there are nilpotent groups of class $n$ that are not $(n-1)$-Engel. An example is the group $G=C_p \wr C_{p}$, that is the wreath product of two cyclic groups of order a prime $p$. One can see that $G$ is nilpotent of class $p$ but not a $(p -1)$-Engel group (see Theorem 6.2 of \cite{Liebeck}).
Furthermore is not difficult to construct $n$-Engel groups that are not nilpotent. Take for example a prime $p$ and consider the group $G=C_p \wr C_{p^{\infty}}$, where $C_{p^{\infty}}$ is the Pr\"{u}fer group. One can prove that $G$ is a $(p+1)$-Engel group but it is not nilpotent.

The main open question in the realm of Engel groups is whether every $n$-Engel group is locally nilpotent. Hence, the Engel condition is more general than nilpotency. Recall that a group $G$ is locally nilpotent if each finitely generated subgroup of $G$ is nilpotent. 

For  $n \leq 4$ every $n$-Engel group is locally nilpotent. Indeed, $1$-Engel groups are abelian. For $n=2$, Levi proved that $G$ is a $2$-Engel group if and only if the normal closure of an arbitrary element is abelian \cite{Levi}. Moreover $2$-Engel groups are nilpotent of class at most $3$. For $n=3$, Heineken proved that every $3$-Engel is locally nilpotent \cite{heineken1961}. Finally, Havas and Vaughan-Lee \cite{Havas2005} proved that every $4$-Engel groups generated by $d$ elements are nilpotent of class $4d$. For $n>4$ the question is still open. 

In some other classes of groups also being $n$-Engel implies locally nilpotency. For example, any residually finite $n$-Engel group is locally nilpotent \cite{wilsonresfinite}.

About Engel groups, as pointed out in the introduction, the question whether every Engel group is locally nilpotent is the analogue of the General Burnside Problem for Engel groups.

It is easy to check that every locally nilpotent group is Engel. The converse is true in some classes of groups like finite groups (Zorn \cite{Zorn}), soluble groups (Gruenberg \cite{Gruenberg}), groups with maximal condition (Baer \cite{Baer}), linear groups (\cite{linearengelgroups}) and some other classes of groups. Hence, in these cases the condition to being Engel is equivalent to local nilpotency, but in general it is much weaker. Indeed, the infinite $p$-group $G = \langle x_{1} \dots, x_{r} \rangle$ constructed by Golod is a counterexample. This group is Engel and it is such that every $(r-1)$-generated subgroup is nilpotent. So far, in the periodic case, Golod's group is the only known example of finitely generated Engel group that is not nilpotent.

\subsection{Identities satisfied by $n$-Engel groups}\label{identityengel}

In this section we present some semigroup identities satisfied by $n$-Engel groups for  $n \leq 4$.
We have the following.
\begin{enumerate}
    \item In every $2$-Engel group $yx^2y=xy^2x$.
    \item In every $3$-Engel group $$xy^2xyx^2y = yx^2yxy^2x \mbox{ and } xy^2xyxyx^2y = yx^2y^2x^2y^2x.$$
    \item In every $4$-Engel group we have
    $$xy^2xyx^2y^2x^2yxy^2xyx^2yxy^2x^2y^2xyx^2y =yx^2yxy^2x^2y^2xyx^2yxy^2xyx^2y^2x^2yxy^2x.$$
\end{enumerate}
For more information one can see \cite{nroot} and \cite{semigroupidentities}.

Notice that, additionally, if a group is Engel and locally nilpotent one can prove the following.

\begin{thm}[\cite{crosbytraustason}]\label{identitytraustason}
There exist positive integers $m = m(n)$ and $r = r(n)$ such that any locally nilpotent $n$-Engel group satisfies
$$
[x^r, x_1, \dots, x_m]=1.
$$
\end{thm}

\subsection{Burnside groups}
In this subsection, we present Burnside groups since it seems that their behaviour is similar to the behaviour of the class of Engel groups. Moreover, in Section \ref{learningburnside} we present a protocol that uses Burnside groups of exponent 3 as platform group.

Let $F_m$ be a free group of rank $m$. The free Burnside group $B(m,n)$ is the group $F_m/F_m^n$, where $F_m^n$ is the group generated by all the $n$-th powers of elements of $F_m$. Therefore $B(m,n)$ is the group in which the identity $x^n=1$ holds and for this reason is the biggest group generated by $m$ elements of exponent $n$. The Burnside Problem asks if a finitely generated group of finite exponent is finite, that is the same as asking whether the free Burnside groups  are finite.

It is easy to prove that for any $m$ the $2$-group $B(m,2)$ is elementary abelian and so finite. Burnside proved that also $B(m,3)$ is finite for any $m$. Levi and van der Waerden in 1993 proved that if $m \geq 3$ the Burnside group $B(m,3)$ is finite and nilpotent of class $3$. One can also prove that since $B(m,3)$ is a group of exponent $3$, then it is a $2$-Engel group because of the following.
\begin{prop}[12.3.5 of \cite{rob}]\label{burnside2engel}
A group of exponent 3 is a 2-Engel group.
\end{prop}

In 1940, Sanov proved that also the group $B(m,4)$ is finite for every $m$. For $n=5$ the problem is still open. Although for some small values of $n$ the Burnside groups are finite, in 1968,  Novikov and Adian proved that, in general, Burnside groups need not to be finite. Indeed they showed that if $m \geq 2$, and $n$ is odd and greater than $4381$, then $B(n,m)$ is infinite. This bound was improved later by Adian which showed that $n$ can be chosen odd and greater than or equal to $665$. With the same construction used to prove that these Burnside groups are infinite,  Novikov and Adian also proved that the word and the conjugacy problems are decidable in Burnside groups, see \cite{Novikov}.

\section{\bf Other properties and problems around Engel groups}\label{otherproperties}

In this section, we survey additional properties that can help us to the study of Engel groups. In particular, we present the growth rate, the Discrete Logarithm Problem for finite $p$-groups, and the degree of nilpotency. 

\subsection{Growth rate}

Let $G$ be a finitely generated group. The growth function $\gamma$ is
\begin{align*}
\gamma: \mathbb{N} & \to \mathbb{R}  \\
n & \mapsto |\{w \in G: l(w) \leq n \}|
\end{align*}
where $l(w)$ denotes the length of $w$. Since  words are used as keys in group-based cryptography, there is a natural relationship between the growth rate of a group and the key space, that is the set of all possible keys. A fast growth rate engenders a large key space, making an exhaustive search of this space intractable. 
Notice that by a deep result of Gromov, if $G$ is a finitely generated group, then $G$ has polynomial growth if and only if $G$ is virtually nilpotent \cite{Gromov}. Moreover Golod Shafarevich groups have exponential growth \cite{growthgsgroups}.

\subsection{The Discrete Logarithm Problem}
%Let $G$ be a group and let $x, y \in G$. The Discrete Logarithm of $y$ to the base $x$ is the least positive integer $a$ for which $x^a=y$.
Let $G$ be a finite group. Given $x, y \in G$, the Discrete Logarithm Problem (in the reminder DLP for short) is to find a positive integer $a$ such that $x^a=y$ (if it exists). Notice that this is usually defined in the setting of cyclic groups because the Discrete Logarithm exists for all elements and all nontrivial bases. 

Notice that the DLP can be generalized to several components as follows. Let $\textbf{x} = (x_1, \dots, x_n)$ be a tuple of elements such that $G = \langle x_1, \dots, x_n \rangle$. Given $y \in G$, the ``generalized'' DLP of $y$ with respect to $\textbf{x}$ is to find $a_i$ such that $y$ can be written uniquely as $\textbf{x}^{\textbf{a}} = (x_1^{a_1}, \dots, x_n^{a_n})=y$ where $0 < a_i < |x_i|$ for every $i$.

For a survey on the topic see, for example, \cite{dlpsurvey1}.

\subsection{Degree of $n$-nilpotency}\label{degreeofnilp}

In this subsection we present the definition of degree of $n$-nilpotency and some related properties. This is a notion that measures how close is a group to being nilpotent of nilpotency class $n$. In particular, this is a way to measure how close Engel groups are to being nilpotent. This is relevant for us since some protocols work better using Engel (or nilpotent) groups (see Section \ref{applications}). In the following we provide two ways to compute the degree of nilpotency of a group.

\subsubsection{Degree of $n$-nilpotency of a finite group}\label{degreeERL}

In \cite{MCS} the degree of $n$-nilpotency ($n \geq 2$) of a finite group $G$ is defined as follows:
\begin{eqnarray*}
d^{(n)}(G) = \frac{|\{(x_1, \cdots, x_{n+1})\in G^{n+1} \mid [x_1, \cdots, x_{n+1}]=1\}|}{|G|^{n+1}}.
\end{eqnarray*}

For ease of notation write $d(G)=d^{(1)}(G)$. Clearly $d(G) = 1$ if and only if the group $G$ is abelian. In \cite{Gustafson} Gustafson proved that if $d(G) > 5/8$, then $G$ is abelian. Moreover if $d(G) > 1/2$, then $G$ is nilpotent \cite{lescot}.

In the following we present two upper bounds for the above definition of $n$-nilpotency degree \cite{msc2}.

\begin{thm}
Let $G$ be a finite group which is not nilpotent of class at most $n$. Then 
$$
d^{(n)}(G) \leq \frac{2^{n+2}-3}{2^{n+2}}.
$$
\end{thm}

\begin{thm}
Let $G$ be a nontrivial finite group with trivial center. Then for every $n \geq 1$,
$$
d^{(n)}(G) \leq \frac{2^{n}-1}{2^{n}}.
$$
\end{thm}

\subsubsection{General definition}
We conclude this section by pointing out that there exists a more general definition of degree of nilpotency which holds also for infinite groups. 
\begin{defn}\cite{degreeofnilp}
Let $G$ be a group (not necessarily finite) generated by a finite set  $X$. The degree of $n$-nilpotency of $G$ with respect to  $X$ is defined as follows:
$$
d_{n,X}(G) = \limsup_{m \rightarrow \infty} \frac{|\{\mathbf{x} \in B_X(m)^{n+1} \mid [x_1, \dots, x_m]=1\}|}{ |B_X(m)|^{n+1}},
$$
where $B_X(m)$ is the ball of radius $n$ centered
at 1 on the Cayley graph of $G$ with respect to $X$.
\end{defn}

Obviously, if the group $G$ is finite then the above definition coincides with the one given in Section \ref{degreeERL}. We finally remark that in \cite{degreeofnilp} it has been proved that the positivity of $d_{n,X}(G)$ does not depend on the generating set $X$.

\section{\bf Algorithmic problems}\label{Sectionalgorithmicproblems}

In this section we survey several decision problems in group theory. In particular, we summarize the status of the complexity of some algorithmic problems in the context of Engel groups. In Section \ref{applications} we present cryptosystems whose security depends on some of these problems.

Notice also that, as pointed out in Section \ref{preliminaries}, every finitely generated $n$-Engel group with $n \leq 4$ is nilpotent. Then every result for nilpotent groups can be applied to $n$-Engel groups for $n \leq 4$.

\subsection{Decision problems}\label{decisionproblems}

From now on we let $G$ be a group given by a presentation $\langle X|R \rangle$ and we understand that when we speak of elements of $G$ these are given as a product of generators in $X^{\pm 1}$.

The following three decision problems were introduced by Dehn in 1911. They are defined as follows.

\begin{itemize}
    \item \emph{Word Problem}: For any $g \in G$, determine if $g$ is the identity element of $G$.
    \item \emph{Conjugacy Problem}: For any $x, y \in G$, determine if $x$ and $y$ are conjugate.
    \item \emph{Isomorphism Problem}: Let $G$ and $G'$ be groups given by finite presentations, determine if $G$ is isomorphic to $G'$. 
\end{itemize}

For polycyclic groups all three of the above problems are decidable (see \cite{conjproblempolyciclic2}, \cite{conjproblempolyciclic}, \cite{conjsep}, \cite{Segalpoly}, and for a survey on the topic see \cite{polycyclic}).

For finitely generated nilpotent groups the word problem is solvable. Also, finitely generated nilpotent groups are linear (\cite{phall}) and for linear groups the word problem is solvable in logspace (\cite{linearlogspace}). Hence, the word problem for finitely generated nilpotent groups is is solvable in (deterministic) logspace.

Furthermore, Blackburn in 1965 proved that the conjugacy problem is decidable as well \cite{conjnilp}.

The isomorphism problem is also known to be solvable for finitely generated nilpotent groups \cite{grunewald}. Notice that, in contrast, the epimorphism problem is undecidable for finitely generated nilpotent groups \cite{Remeslennikov}.
We record that the \emph{epimorphism problem} asks whether, given two finite presentations $R_1$ and $R_2$ of groups, there exists an algorithm which determines if the group with presentation $R_1$ is a homomorphic image of the group defined by the presentation $R_2$.

\subsubsection{Power decision problem}
Let $G$ be a finitely generated group. Given $x, y \in G$, determine whether there exists $n \in \mathbb{Z}$ such that $y = x^n$. Notice that this is equivalent to deciding whether $y \in \langle x \rangle$, and thus is the analogue of the Discrete Logarithm Problem (see Section \ref{DLP}). A generalization of the power problem is the subgroup membership problem that can be stated as follows.

\subsubsection{Subgroup membership decision problem}
Let $G$ be a group and $H \leq G$. The subgroup membership decision problem (also known as the \emph{generalized word problem}) asks for any $g \in  G$ if $g \in H$. The subgroup membership problem is decidable for finitely generated nilpotent groups \cite{submembnilp}.
This problem can be generalized to the \emph{rational subset membership problem}. Recall that given a group $G$, the class of rational subsets of $G$ is the smallest class that contains all finite subsets of $G$ and that it is closed with respect to
union, product and taking the free monoid generated by a set (i.e.\ using the Kleene star operation). The rational subset membership problem asks whether given a rational subset $H$ of $G$ and an element $g \in G$, whether $g \in H$. In this case, nilpotent groups have undecidable rational subgroup membership problem \cite{rationalsubset}.

\subsubsection{The endomorphism decision problem}
Let $F$ be a free group and $x \in F$. The endomorphism decision problem asks whether there exists an algorithm that given $y \in F$ decides if there exists an endomorphism $\phi$ of $F$ sending $x$ to $y$. In nilpotent groups the endomorphism problem is undecidable \cite{endomnilp}.

\subsubsection{The $n$-th root decision problem}

For $n \in \mathbb{N}$, an element $x$ of a group $G$ is a $n$-th root of an element $a$ if the equation $x^n=a$ has solution.  Given an element $g \in G$, the $n$-th root decision problem asks to determine if $g$ has any $n$-th root. Note that if the root of every element of $G$ belongs to $G$, the group is said to be \emph{complete}. 

We remark that if we take a group $H$ under addition, this is the problem of finding $n$-th roots is equivalent to saying that for any $n \in \mathbb{N}$ and every element $h \in H$, the equation $nx = h$ has at least one solution in $H$.

For some results about finding $n$-th roots in nilpotent groups see \cite{nroot}. It seems that finding square roots in Engel groups is a difficult problem (see Section \ref{2eng} and Section \ref{4eng}).

\subsubsection{Geodesic length decision problem}

Let $G = \langle X \rangle$ be a group generated by a set $X$. Denote  with $|w|$ the length of a word in the alphabet $X^{\pm 1}$, and with $\rho$ the canonical epimorphism of $F(X)$ onto $G$, where $F(X)$ denotes the free group over the set $X$. The geodesic length $l_X(g)$ of an element $g \in G$ with respect to $X$ is 
$$l_X(g) = \min\{|w| \mid w \in F(X), \ \rho(w) = g\}.
$$
We say that a word $w \in F(X)$ is  geodesic if $|w| = l_X(\rho(w))$. Given a word $w \in F(X)$, the geodesic length decision problem consists in, given a word $w \in F(X)$,  to find $l_X(\rho(w))$. 

In \cite{Myasnikovbook} it is observed that for nilpotent groups this is a hard problem.

\subsection{Search problems}\label{searchproblems}
In this section we survey some search  problems. Several protocols of non-commutative cryptography, like Ko-Lee and Anshel-Anshel-Goldfeld, are based in part on the conjugacy search problem. For this reason here we present, among others, the conjugacy search problem and some of its variations. 

As before, we let $G$ be a group given by the presentation $\langle X|R \rangle$ and we agree that when we refer to an element of $G$ this is given as a product of generators in $X^{\pm 1}$.

\subsubsection{Word search problem}
 The word search problem is: given a finitely presented group $G$ and an element $g = 1$  in $G$ find a presentation of $g$ as a product of conjugates of defining relators and their inverses.
 
\subsubsection{Conjugacy search problem} 
Let $G$ be a group and $a_1, \dots a_n, b_1, \dots b_n \in G$ where $a_i$ is conjugate to $b_i$ for $i=1, \dots, n$. The multiple conjugacy search problem  is to find $c \in G$ such that $a_i^c = b_i$ for $1 \leq i \leq n$. If $n=1$ this reduces to the conjugacy search problem. 

In some cases, for example for polycyclic groups  \cite{Eickconj}, the multiple conjugacy search problem reduces to the solvability of single (independent) conjugacy search problem. In general, for finitely generated polyciclic groups the conjugacy search problem can be solved  by recursively enumerating the conjugates of the elements taken in consideration \cite{searchvladimir}. 

\subsubsection{Power conjugacy search problem}

The power conjugacy search problem asks to find for some $x, y \in G$, an element  $g \in G$, and $n \in \mathbb{N}$ such that  $x^n =y^g$. Note that for $n = 1$ this reduces to the standard conjugacy search problem, whereas if $g$ is the identity in the group this reduces to the power (search) problem.

\section{\bf Applications to cryptography}\label{applications}

In this section we present some applications of Engel groups to cryptography. We first survey some cryptosystems based on $n$-Engel groups such as a public key, a digital signature, two secret sharing schemes and multilinear maps. We present the protocol used, the ideal platform group, and some security assumptions. We conclude this section presenting the Discrete Logarithm Problem (DLP for short) in finite $p$-groups (that are nilpotent groups, so in particular Engel groups) and providing the actual status of quantum algorithms.

\subsection{Multilinear maps}\label{multilinearmaps}

In this subsection we first define multilinear maps and then we provide some applications. For more information, see \cite{multilinear}.

\begin{defn}
Let $n$ be a positive integer and $G$ an arbitrary group. A map $e: G^n  \to G$ is said to be a $n$-linear map (or a \emph{multilinear map}) if for any $g_1, \dots, g_n, \in G$ and any $a_1, \dots, a_n \in \mathbb{Z}$ we have 
    $$
    e(g_1^{a_1}, \dots, g_n^{a_n}) = e(g_1, \dots, g_n)^{a_1 \cdots a_n}.
    $$
Moreover, we say that the map $e$ is non-degenerate if there exists $g\in G$ such that $e(g, \dots g) \neq 1$.

\end{defn}

Let now $G$ be a nilpotent group of class $n > 1$ and  $g_1, \dots, g_n \in G$. One can easily prove by induction on $n$ that for any $a_1, \dots, a_n \in \mathbb{Z}$ the following identity holds:
\begin{equation}\label{identitynilpotent}
   [g_1^{a_1}, \dots, g_n^{a_n}] = [g_1, \dots, g_n]^{\prod_{i=1}^n a_i}.  
\end{equation}

Hence if $G$ is nilpotent, the map $e$ 
\begin{align*}
 e: G^{n} & \to G \\
 (g_1, \dots, g_{n}) & \mapsto [g_1, \dots, g_{n}]
\end{align*}
is a multilinear map. In addition, if we fix $x \in G$, we can construct another multilinear map $f$ given by
\begin{align*}
 f: G^{(n-1)} & \to G \\
 (g_1, \dots, g_{n-1}) & \mapsto [x, g_1, \dots, g_{n-1}].
\end{align*}
Notice that we want the multilinear map to be non-degenerate. This is the case if we assume that $G$ is not $(n-1)$-Engel. Indeed, in this condition, there exists $g \in G$ such that $[x,_{(n-1)} g] \neq 1$, which implies that $f$ is non-degenerate.

\subsubsection{Protocol}
Let $G$ be a nilpotent group of class $n + 1$ which is not $n$-Engel  ($n \geq 1$). Then there exist elements $x, g \in G$ such that $[x,_n g]\neq 1$. The group $G$ is public. We have $n + 1$ users $A_1, \dots A_{n+1}$ that wish to agree on a shared secret key. Each user $A_j$ selects a private integer $a_j \neq 0$, computes $g^{a_j}$, and sends it to the other users. Then we are in the following situation:
\begin{itemize}
    \item The user $A_1$ computes $[x^{a_1},g^{a_2}, \dots g^{a_{n+1}}]$.
    \item For $j=2, \dots, n$, the user $A_j$ computes $[x^{a_j}, g^{a_1}, \dots, g^{a_{j-1}}, g^{a_{j+1}}, \dots, g^{a_{n+1}}]$. 
    \item The user $A_{n+1}$ computes $[x^{a_{n+1}},g^{a_1}, \dots g^{a_{n}}]$.
\end{itemize}

Since \eqref{identitynilpotent}  holds in all nilpotent groups, all elements computed by the users are equal to $k=[x,_n g]^{\prod_{j=1}^{n+1} a_j}$. This is the shared key.

\subsection{Two cryptosystems based on $n$-th root problem in $n$-Engel groups}\label{cryptoengel}

In this section we present two cryptosystems in $n$-Engel groups proposed in \cite{nroot}.

\subsubsection{A public key based on $2$-Engel groups}\label{2eng}

A satellite generates some data from a 2-Engel group. Alice and Bob choose two elements $x$ and $y$ as their secret keys, respectively. Then Alice sends
$x^2$ to Bob and Bob sends $y^2$ to Alice. As pointed out in Section \ref{identityengel},  in any $2$-Engel group it holds the semigroup identity $xy^2x = yx^2y$. Then Alice and Bob agree on a key. One could extend this scheme to other $n$-Engel groups, since there are similar relations in them as well (see, for example Equation \ref{identitytraustason}).

\subsubsection{A digital signature based on $4$-Engel groups}\label{4eng}

Consider a $4$-Engel group, which is nilpotent and satisfies the following semigroup law (see the last part of Section \ref{Engelgroups} and Section \ref{identityengel}):
$$
xy^2xyx^2y^2x^2yxy^2xyx^2yxy^2x^2y^2xyx^2y =yx^2yxy^2x^2y^2xyx^2yxy^2xyx^2y^2x^2yxy^2x.
$$
The idea to make a digital signature is as follows. Suppose $x$ and $y$ are private and $x^2$, $y^2$, $xy^2x$, and $xy^2x^2y^2x$ are public.
The public key is $x^2$ and the signature is $xy^2x$ and $xy^2x^2y^2x$. The verifier knows $y$, so he has to verify both the semigroup identity. 

\subsubsection{Platfrom group and security} 
The underlying group could be combined by a finite group in which square root is hard, for example $\mathbb{Z}^*_{pq}$, with $p, q$ prime numbers.

The security of these digital signatures lies on the fact that the complexity of
finding square root in $2$-Engel and $4$-Engel groups. It seems that this is a hard problem \cite{computationnilpotent}.

\subsection{Secret sharing based on the word problem in Engel groups}

In this section we present two secret sharing schemes based on the word problem proposed by Habeeb-Kahrobaei-Shpilrain \cite{secretsharing}. Both of them are based on the efficiency of the word problem and we will use them with ($n$-)Engel nilpotent groups.

We recall that a $(t,n)$-threshold secret sharing is a scheme in which a secret is distributed among $n$ participants in a way that the secret can be recovered only if at least $t$ of them combine their shares.  We present two secret sharing schemes. The first one (Scheme 1) is a $(n,n)$-threshold scheme and so the participants must get together to recover the secret. The second one (Scheme 2) is a $(t,n)$-threshold scheme that is a combination of Shamir’s scheme \cite{shamir} and the group-theoretic scheme proposed in Section 3 of \cite{secretsharing} and presented below. In both of them we denote a generic participant of the protocol with $P_i$, with $1 \leq i \leq n$. We also suppose that the dealer and the participants at the beginning are able to communicate over secure channels and then they communicate over open channels.

\subsubsection{Scheme 1}
The following is an $(n,n)$-threshold scheme.
The dealer:
\begin{itemize}
\item Distributes a $k$-column  of bits   $C = (c_1, c_2, \dots, c_k)^{\textsf{T}}$.
\item Distributes $C$ among $n$ participants  in  such  a  way  that the column can be reconstructed only if all participants combine their information.
\end{itemize}
A set $X=\{x_1, \dots, x_m\}$ of generators is public.  Then the protocol is as follows.

\begin{enumerate}
    \item  The dealer uses a secure channel to assign to each participant $P_j$ a set of words $R_j$ in the alphabet $X^{\pm 1}$ such that each group $G_j = \langle x_1, \dots, x_m \mid R_j \rangle$ has word problem solvable in polynomial time.
    \item  The dealer splits the column $C$ in $\sum_{j=1}^n C_j \mod 2$. These are secret shares to be distributed to the $n$ participants.
\item The dealer distributes words $w_{1j}, \dots w_{kj}$ in the generators $x_1, \dots x_m$ over an open channel to each participant $P_j$, with $1 \leq j \leq n$. The choice of the words is such that $w_{ij} \neq 1$ in $G_j$ if $c_{ij} =0$ and $w_{ij} =1$ in $G_j$ if $c_{ij} =1$, where $c_{ij}$ is the $i$-th entry of $C_j$.
\item Each participant $P_j$ check for any $i$ if $w_{ij}=1$ in the group $G_j$ or not.
\item Now each participant constructs $C_j = (c_{1j} , c_{2j}, \dots, c_{kj})^{\textsf{T}}$ of $0$'s and $1$'s ($c_{ij}= 1$ if $w_{ij}=1$ in $G_j$, and $0$ otherwise).
\item The secret is built  by putting together the vector sum $\sum_{j=1}^n C_j \mod 2$.
\end{enumerate}

\subsubsection{Scheme 2}

In this case we present a $(t,n)$-threshold scheme in which the secret is an element $x \in \mathbb{Z}_p$, for $p$ a prime number. The dealer:
\begin{itemize}
    \item Chooses a polynomial $f$ of degree $t-1$ such that $f(0) = x$.
    \item Determines integers $y_i = f(i) \mod p$ for $1 \leq i \leq n$.
    \item Distributes every $y_i$ to the correspondent participant $P_i$.
\end{itemize}
We now consider a  set $X=\{x_1, \dots, x_m\}$ of generators, which is public. In addition, we assume that the integer $x$ and every $y_i$ can be written as $k$-bit columns.
Now we are ready to present the scheme, that reads as follows.
\begin{enumerate}
    \item The dealer distributes over a secure channel to each participant $P_j$ a set of relators $R_j$ such that each group $G_j = \langle x_1, \dots, x_m \mid R_j \rangle$ has efficiently solvable word problem.
    \item The dealer then distributes over open channels $k$-columns of the form
    $b_j = (b_{1j}, b_{2j}, \dots, b_{kj})^{\textsf{T}}$, (where $1 \leq j \leq n$) of words in $x_1, \dots, x_m$ to each participant. The $b_{ij}$ are chosen such that if we replace them by bits, the resulting bit column represents the integer $y_j$. Note that we use ``1'' if $b_{ij} = 1$ in the
group $G_j$ and ``0'' otherwise.
\item For any word $b_{ij}$, the participant $P_j$ checks whether or not $b_{ij} = 1$ in the group $G_j$. Then $P_j$ obtains a binary representation of the number $y_j$, and therefore $P_j$ recovers $y_j$.
\item Each participant now has a point $f(i) = y_i$ of the polynomial. Using polynomial interpolation, any $t$ participants can now recover the polynomial $f$. Whence they obtain the secret $x = f(0)$.
\end{enumerate}
If $t \geq 3$, Step 4 can be changed in a way that the participants do not have to reveal their individual shares to each other if they do not want to. For more details about how to arrange this, an interested reader can see \cite{secretsharing}.

\subsubsection{Platform group}
For Scheme 1 and Scheme 2, we propose finitely presented nilpotent Engel groups that, being nilpotent, have efficiently solvable word problem (see the discussion in Section \ref{decisionproblems}). For example, one can use $n$-Engel groups with $n=2,3,4$ or an Engel group that is linear, solvable or finite (and hence nilpotent). 

\subsection{A key protocol based on semidirect products of groups}
In this section we present a key protocol based on semidirect products of groups (or semigroups). This protocol can be based on any group, but we present an application using free nilpotent $p$-groups where $p$ is a sufficiently large prime $p$ proposed by Kahrobaei and Shpilrain in \cite{semidirect}. Note that there are several protocols based on semidirect products of groups, such as a public key exchange \cite{semidirect2}.

\subsubsection{Protocol}
Let $G$ be a group (or a semigroup) and choose public $g \in G$ and $\varphi \in \Aut(G)$ (or $\End(G)$). Alice chooses a private $m \in \mathbb{N}$ while Bob chooses a private $n \in \mathbb{N}$.
Alice and Bob are going to work with elements of the form $(g, \varphi^r)$, where $g \in G$ and $r \in \mathbb{N}$. Recall that the multiplication of two elements of this form is as follows:
$$
(g, \varphi^r)(h, \varphi^s)=(\varphi^s(g)h, \varphi^{r+s}).
$$
Now we present a key exchange protocol similar to the Diffie-Hellman key exchange.
\begin{itemize}
    \item Alice computes 
    $$(g, \varphi)^m = (\varphi^{m-1}(g) \cdots \varphi^2(g)\varphi(g)g, \varphi^m)
    $$
    and sends the first component $a=\phi^{m-1}(g) \cdots \varphi^2(g)\varphi(g)g$ of the pair to Bob.
     \item Bob computes 
    $$(g, \varphi)^n = (\varphi^{n-1}(g) \cdots \varphi^2(g)\varphi(g)g, \varphi^n)
    $$
    and sends the first component $b=\varphi^{n-1}(g) \cdots \varphi^2(g)\varphi(g)g$ of the pair to Alice.
    \item Alice computes 
    $$(b,x)(a, \varphi^m)=(\varphi^{m}(b)a, x\varphi^m).
    $$
    Her key is $K_A=\varphi^{m}(b)a$.
     \item Bob computes 
    $$(a,y)(b, \varphi^n)=(\varphi^{n}(a)b, y\varphi^n).
    $$
    His key is $K_B=\varphi^{n}(a)b$.
\end{itemize}
The shared secret key is $K=K_A=K_B=$ since 
    $$
    (b,x)(a, \varphi^m)=(a,y)(b, \varphi^n)=(g,\varphi)^{m+n}.
    $$

\subsubsection{Platform groups}
We denote with $\gamma_c(G)$ the normal subgroup of $G$ generated by all the elements of the form $[y_1, \dots, y_c]$.
Consider now $F_m$ the free group on $x_1, \dots, x_m$.  The factor group $F_m/\gamma_{c+1}(F_m)$ is the free nilpotent group of nilpotency class $c$. 

The group suggested as platform group for the protocol described above is 
$$
G=F_m/F_m^{p^2}\gamma_{c+1}(F_m)
$$
that is a nilpotent $p$-group and hence a finite group whose order depends on $c, m$ and $p$. The suggested values of $c$ and $m$ are small numbers, for example $c=2$ or $c=3$ to make the computation efficient. Conversely, the value of $p$ should be sufficiently large to make it secure to a linear algebra attack.

\subsection{Burnside groups of exponent $3$ and the learning homomorphism problem}\label{learningburnside}

The learning with errors (LWE) problem was introduced by Regev in \cite{Regev}, and has become one of the most known problems in lattice-based cryptography. It has been used to construct several cryptosystems, and it is believed to be hard even for quantum computers. 

Informally speaking, the problem of learning with error is to deduce a particular function by sampling the input/output behavior if some of the outputs are incorrect. For more information one can see \cite{learningproblem}. 

In the following, we define the Learning Homomorphim Problem with Noise (LHN, for short) that is a generalization of the LWE. We will denote the set of homomorphisms from a group $G$ to a group $H$ with $\Hom(G, H)$.

\begin{defn}
Let $\varphi \in \Hom(G, H)$, where  $G$ and $H$ are finitely generated groups, and let $g_1, \dots g_m \in G$. Also, let $\alpha$ and $\beta$ be probability distributions over $G$ and $H$, respectively. Let $\Psi$ be the probability distribution over $G \times H$ which assigns to each tuple of the form $(g, \varphi(g)h)$ the probability $\Prob(g\sim \alpha)$ and $\Prob(h\sim \beta)$ and the rest of tuples from $G\times H$ are assigned a probability of $0$.

We say that an algorithm solves the LHN with noise if for any $\varphi: G \to H$, the algorithm is able to learn $\varphi$ given a set of samples from the distribution $\Psi$ with high probability. For the purpose of this paper we omit the definition of ``learning with high probability given a set of samples''. For a detailed reference on this, see \cite{learningproblem} and \cite{learningproblem2}.
\end{defn}

\subsubsection{Protocol, platform groups and security assumptions}

In \cite{learningproblem2} the authors employed for the first time properties of Burnside groups in cryptography.  In particular, they introduced the Learning \emph{Burnside} Homomorphisms with Noise (B\--LHN), which uses for computational purposes only surjective homomorphisms in contrast to the general definition of LHN \cite{learningproblem}. 

To avoid technicality, we omit the description of the protocol used for this problem. We address the reader to pp.\ 10-11 of \cite{learningproblem}.

We want to underline that as platform group for the B-LHN problem, the authors employed the Burnside groups of exponent $3$, that are, in particular, $2$-Engel groups. The security of this protocol is based on the computational hardness of B-LHN (see Theorem 2 of \cite{learningproblem2}).

\subsection{The Discrete Logarithm Problem in finite $p$-groups}\label{DLP}

Sutherland studied the DLP in some finite abelian $p$-groups \cite{shuterland}, and in a series of papers by Mahalanobis, the DLP has been studied for finite $p$-groups of nilpotency class 2 \cite{DLP2}. Solving the DLP in finite $p$-groups of larger class is an interesting question. One, for example, can consider  semidirect product of cyclic $p$-groups of well-defined orders, to make a nilpotent group and then computing the Discrete Logarithm in each factor.

\subsection{Quantum algorithms}\label{quantum}
We conclude this section of applications by presenting the status of quantum algorithms in cryptography. In 2015 the National Security Agency announced plans to replace all deployed cryptographic protocols with quantum secure protocols. A quantum computer is able to perform integer factorization and solve the DLP in finite cyclic groups in polynomial time.

In \cite{batty} it has been explored the application of quantum algorithms to group theory. In particular, from a group theoretic point of view, Shor's algorithm can be seen as the hidden subgroup problem in finite cyclic groups. We recall that a subgroup $H$ of a group $G$ is \emph{hidden} by a function $f$ from $G$ to a set $X$ if it is constant over all cosets of $H$, and takes distinct values on distinct cosets. In other words, for any $g_1, g_2 \in G$, $f(g_1) = f(g_2)$ if and only if $g_1H = g_2H$.
Given a hidden subgroup $H$, the \emph{hidden subgroup problem} asks to find a generating set for $H$ using information from evaluations of $f$ via an oracle (for a survey on the topic, see \cite{Delahoran}). In \cite{hspsemidirect} the complexity of quantum algorithms for the HSP for certain semidirect product of certain finite $p$-groups, for $p$ a prime has been considered.

\section{\bf Conclusion and open questions}\label{conclusion}

In this paper we have presented the current state of Engel group-based cryptography. We have begun with a study of algorithmic properties of Engel groups, and we have also seen that there are a variety of key exchanges, digital signature systems, and secret sharing schemes for which an Engel (or nilpotent) group is an appropriate choice of platform group.
As we have seen throughout the paper, if on the one hand there has been some results about Engel groups and their attendant cryptosystems over the last decade, on the other hand the majority of computational complexity and algorithmic questions remain unanswered. 

We collect some open problems below with the hope of stimulating interest in their solutions.

We start by pointing out that since Engel groups are a generalization of nilpotent groups, in the case in which the algorithmic problems presented in in Section \ref{Sectionalgorithmicproblems} are decidable for nilpotent groups, we ask whether the same holds specially if the group is $n$-Engel. Also, since some problems are undecidable for nilpotent groups (take, for example, the rational subgroup membership problem and the endomorphism problem) we ask:
\begin{enumerate}
    \item Is the rational subgroup membership problem undecidable for Engel groups?
    \item Is the endomorphism problem undecidable for Engel groups? 
\end{enumerate}  

In addition, we propose the following questions related also to the complexity of some algorithmic problems presented in Section \ref{Sectionalgorithmicproblems}. 
\begin{enumerate}[resume]
\item Determine if the following algorithmic problems are decidable and if yes, find the complexity:
\begin{itemize}
\item The word problem.
\item The (power) conjugacy problem.
\item The geodesic length problem.
\end{itemize}
\item What about the search problems proposed in Section \ref{searchproblems}?
For example, if $G$ is an Engel group, and $a_1, \dots a_n, b_1, \dots b_n \in G$ where $a_i$ is conjugate to $b_i$ for $i=1, \dots, n$. Can we find $c \in G$ such that $a_i^c = b_i$ (for $1 \leq i \leq n$)? As a starting point, one can study the single conjugacy search problem (i.e.\ the case $n=1$). 
\end{enumerate}

In Section \ref{applications} we present several ideas of Engel-group based cryptography.

In particular, in Section \ref{multilinearmaps}, we present a protocol based on multilinear maps. It is an open question the ideal platform of this scheme. Thus we ask the following.
\begin{enumerate}[resume]
    \item What is an ideal platform group for the protocol proposed in Section \ref{multilinearmaps}?
\end{enumerate}

In Section \ref{cryptoengel} we present some specific cryptosystems based on $n$-Engel groups. It seems they are secure because they are based on the difficulty to find square roots in $n$-Engel groups. Moreover, no platform groups are suggested. Whence, we propose the following.
\begin{enumerate}[resume]
\item What is an ideal platform group for the protocol proposed in Section \ref{2eng} and Section \ref{4eng}, respectively?
\item  What is the complexity of the square root problem in $n$-Engel groups? And for the $n$-root problem?
\end{enumerate}

In view of Section \ref{quantum}, we conclude the paper with the following questions.
\begin{enumerate}[resume]
\item Is the HSP solvable for Engel groups?
\item Are there some Engel group-based cryptosystems resistant to quantum algorithms? Note that this could be done either by analyzing the HSP for Engel groups or showing the underlying security problem is NP-complete or NP-hard.
\end{enumerate}

\vskip 0.4 true cm

%------------------------------------------------------------------------------------%

%-----------------------------------------------------------------------------
%-----------------------------------------------------------------------------

\bigskip

\section*{Acknowledgements}
The authors would like to thank G. Fern\'andez-Alcober, A. Garreta, and A. Tortora for helpful comments and suggestions. The first author would like to thank A. Tortora and M. Tota for an invitation to the University of Salerno where the ideas of this research were initiated. The second author would like to thank the Department of Computer Science at the University of York for its hospitality while this paper was being written.

\bigskip
\bigskip

\bibliography{bib}

\end{document}